\documentclass[11pt,oneside,english,reqno]{amsart}
\usepackage[T1]{fontenc}
\usepackage[latin9]{inputenc}
\usepackage{color}
\usepackage{babel}
\usepackage{amstext}
\usepackage{amsthm}
\usepackage{setspace}

\makeatletter
\theoremstyle{plain}
\newtheorem{thm}{\protect\theoremname}
  \newtheorem*{thm*}{\protect\theoremname}
  \theoremstyle{plain}
  
\newtheorem*{conj}{Conjecture}
\newtheorem{cor}{Corollary}
\newtheorem{claim}{Claim}
\theoremstyle{remark}

\usepackage{graphicx}
\usepackage{color}

\definecolor{myblue}{rgb}{0.09,0.32,0.44} 

\theoremstyle{remark}
\newtheorem*{qst*}{Question}
\newtheorem*{rmrks*}{Remarks}

\newlength{\tempindent}

\newcommand{\NN}{\mathbb{N}}
\newcommand{\ZZ}{\mathbb{Z}}

\newcommand{\lazyenum}{
\setlength{\tempindent}{\parindent}
\begin{enumerate}[leftmargin=0cm,itemindent=0.7cm,labelwidth=\itemindent,labelsep=0cm,align=left,label=\arabic*)]
\setlength{\parskip}{\smallskipamount}
\setlength{\parindent}{\tempindent}
}

\makeatletter
\renewcommand{\andify}{%
  \nxandlist{\unskip, }{\unskip{} \@@and~}{\unskip{} \@@and~}}
\def\author@andify{%
  \nxandlist {\unskip ,\penalty-1 \space\ignorespaces}%
    {\unskip {} \@@and~}%
    {\unskip \penalty-2 \space \@@and~}%
}
\let\@wraptoccontribs\wraptoccontribs
\makeatother

\makeatother

\providecommand{\lemmaname}{Lemma}
\providecommand{\theoremname}{Theorem}
\providecommand{\theoremname}{Theorem}

\begin{document}

\begin{abstract}
We prove that any Cayley graph $G$ with degree $d$ polynomial growth does not satisfy $\{f(n)\}$-containment for any $f=o(n^{d-2})$. This settles the asymptotic behaviour of the firefighter problem on such graphs as it was known that $Cn^{d-2}$ firefighters are enough, answering and strengthening a conjecture of Develin and Hartke. We also prove that intermediate growth Cayley graphs do not satisfy polynomial containment, and give explicit lower bounds depending on the growth rate of the group. These bounds can be further improved when more geometric information is available, such as for Grigorchuk's group.
\end{abstract}

\title[firefighter lower bounds]{The firefighter problem on polynomial and intermediate growth groups}

\author{Gideon Amir}
\address{Bar-Ilan University, Ramat Gan 52900, Israel.}
\email{gidi.amir@gmail.com}

\author{Rangel Baldasso}
\address{Bar-Ilan University, Ramat Gan 52900, Israel.}
\email{rangel.bal@gmail.com}

\author{Gady Kozma}
\address{The Weizmann Institute of
  Science, Rehovot 76100, Israel}
\email{gady.kozma@weizmann.ac.il}

\maketitle

\section{Introduction}

 Let $G$ be a graph and let $\{f(n)\}$ be a sequence of integers;  an initial fire starts at a finite set of vertices; at each time interval $n\geq 1$, at most $f(n)$ vertices which are not on fire become protected, and then the fire spreads to all unprotected neighbors of vertices on fire; once a vertex is protected or is on fire, it remains so for all time intervals. The graph $G$ has the \emph{$\{f(n)\}$-containment property} if every initial fire admits a strategy consisting of  protecting at most $f(n)$ vertices at the $n^{th}$ time interval so that the set of vertices on fire is eventually constant. We say that the graph $G$ has the \emph{$O(n^d)$-containment property} if there is a  constant $c\geq 1$ such that $G$ has {the $\{cn^d\}$-containment property}. Understanding the relation between $G$ and its containment functions is considered an asymptotic version of the firefighting game introduced by Hartnell in \cite{H}.

\subsection{Background}
 In \cite{DMpT} Dyer, Martinez-Pedroza and Thorne prove the following upper bound:
 \begin{thm*} [\cite{DMpT} Theorem 2.3]
Let $G$ be a connected graph with polynomial growth of degree at most $d$. Then $G$ satisfies the $O(n^{d-2})$--containment  property.
\end{thm*}
This is a generalization of a theorem by Develin and Hartke~\cite{DH}, which showed this holds for $\ZZ^d$. They also prove the following more general upper bound:
\begin{thm*} [\cite{DMpT} Theorem 2.8]
Let $G$ be a locally finite connected graph, let $x_0$ be some vertex in $G$ and let $v(n)= |B_n(x_0)|$ be its rooted growth function. Let $v''(n)$ denote its discrete 2nd derivative $v''(n):=v(n)-2v(n-1)+v(n-2)$. If $v''(n)$ is non-negative and non-decreasing then $G$ has the $\{f_n\}$-containment property with $f_n = 3 v''(2n)$.
\end{thm*}

Essentially the containment strategy behind the above theorems consists of constructing a large (far enough) sphere around the initial fire.
The problem of finding good (possibly matching) lower bounds remained quite open, and only partial results are known.

A first observation is that for general graphs, growth does not capture the correct containment. In fact there are bounded degree graphs of exponential growth where 1 firefighter is enough to contain any initial function. One such graph is the Canopy tree (start with $\NN$ and to each vertex $n\in\NN$ connect a rooted binary tree of depth $n$).

However, one can hope that when the graph is symmetric enough, or specifically on Cayley graphs, one would have a stronger relation  between growth and containment.

Develin and Hartke~\cite{DH} conjectured that for $\ZZ^d$, the following converse holds:
\begin{conj}\cite[Conjecture 9]{DH}
Suppose that $f(n)=o(n^{d-2})$. Then there exists some outbreak on $\ZZ^d$ which cannot be contained by deploying $f (n)$ firefighters at time $n$.
\end{conj}

Partial results were proved by Dyer, Martinez-Pedroza and Thorne  in \cite{DMpT}. To that end they studied a spherical version of isoperimetry, looking at the number of neighbours of subsets of $S_n$ inside $S_{n+1}$. They proved that if $G$ satisfies a "spherical isoperimetric inequality" in the sense that  every   $A\subset S_n$ satisfies $\frac{|\partial A \cap S_{n+1}|}{|A|}\geq\frac{|S_{n+1}|}{|S_n|}\geq 1$, then $G$ does not satisfy $f(n)$-containment for any $f$ for which $\sum \frac{f(n)}{|S_n|}<\infty$. This allows to deduce that $\ZZ^d$ does not satisfy $\{cn^{d-3}\}$-containment, but fails to fully resolve the conjecture. One should also note that many groups do not satisfy this spherical isoperimtric inequality (e.g. the lamplighter group $\ZZ_2 \wr \ZZ$.)

Note that by results of \cite{DMpT} containment is a quasi-isometry invariant, in the following sense: define $g\preceq f$ if there exists $C>0$ s.t. $g(x)\leq Cf(Cx)$ and $g\simeq f$ if $g\preceq f$ and $f\preceq g$. Theorem 4.4. of \cite{DMpT} states that if $G_1$ is quasi-isometric to $G_2$, and $G_1$ satisfies the $\{f(n)\}$-containment property, then $G_2$ satisfies the $\{g(n)\}$-containment property for some $g\simeq f$. It is well-known that for groups, growth and isoperimetric profiles are also invariant under quasi-isometries and, in particular, their equivalence classes do not depend on the choice of finite symmetric generating set.
\subsection{Results}

We are now ready to state our main Theorem, which settles the Develin-Hartke conjecture. In fact, we show that the conjecture holds not only for $\ZZ^d$ but for any Cayley graph of polynomial growth.

\begin{thm}\label{thm:main}
Let $G$ be a Cayley graph satisfying $v(n)\succeq n^d$ for some $d\geq 2$. Then $G$ does not satisfy $\{f(n)\}$-containment for any $f=o(n^{d-2})$
\end{thm}

Our method uses isoperimetry to get lower bounds. By known connections between growth and isoperimetry, we are able to translate this into lower bounds depending on the growth for any growth rate. For exponential growth groups it is not optimal, and in fact does not even give an exponential lower bound. (Whereas by the results of \cite{L}, who uses very different methods, the critical threshold for exponential containment coincides with the growth rate, see also earlier results for elementary amenable groups in \cite{Mp}). However it does allow us to get a superpolynomial lower bound on any group of intermediate growth, answering a well known open question (See e.g. \cite{L} and Question $12$ in \cite{DMpT}).
\begin{cor}\label{c:inter}
Let $G$ be the Cayley graph of some intermediate growth group. Then $G$ does not satisfy polynomial containment of any degree.
\end{cor}

Currently, in all examples of groups of intermediate growth for which concrete lower bounds on the growth are known, these groups satisfy a stretched exponential lower bound of the form $|B_n(id)|\geq Ce^{n^{\alpha}}$ for some $0<\alpha<1$ (in fact, the famous gap-conjecture asserts that all groups of intermediate growth satisfy such a lower bound for some fixed $\alpha_0>0$). Given such a lower bound we can strengthen our result:

\begin{thm}\label{t:strectched}
Let $G$ be a Cayley graph satisfying $v(n)\succeq e^{n^{\alpha}}$ for some $0<\alpha<1$.  Fix any $\beta<\frac{\alpha}{\alpha+1}$, and let $f(n)$ be some function satisfying $f(n)=o(e^{n^\beta})$. Then $G$ does not satisfy $\{f(n)\}$-containment.
\end{thm}

In \cite{EZ}, it is proved an isoperimetric inequality for Grigorchuk groups that improves upon the one attained directly form its growth. We can apply these bounds to get an improved lower bound on containment for Grigorchuk groups. (See \cite{EZ} for exact definitions).
\begin{claim}
Let $G_\omega$ be a Gricorchuk group where $\omega\in \{0,1,2\}^\NN$ is not eventually constant, then $G_{\omega}$ does not satisfy $\{f(n)\}$-containment for any $f(n)=o(e^{n^{1/2}})$.
\end{claim}

\begin{qst*}
Is the above bound tight? That is, does Grigorchuk's group satisfy $\{f(n)\}$-containment for some $f(n)=O(e^{n^{1/2}})$? For some $f(n)=O(e^{n^{1/2+\epsilon}})$ for arbitrary $\epsilon>0$?
\end{qst*}

\section{Proofs}
We begin with some notation:
let $G$ be some finitely generated group and fix some finite symmetric generating set $S$. We will identify $G$ with its Cayley graph w.r.t. $S$. Let $v(n)=|B_n(id)|$ be the volume growth of the group $G$, and let $\Phi$ denote the isoperimetric profile of $G$, that is $\Phi(k)=\inf_{K\subset V, |K|=k}|\partial K|$ (where $\partial K$ denotes the outer boundary of $K$).

Let $K_n$ denote the set of burning vertices at time $n$, and let $k_n:=|K_n|$ denote the number of burning vertices at time $n$.
Let $g(n)$ denote the number of vertices protected up to time $n$. Then $g(n)\leq \sum_{i=1}^n f(i)$.

We now move to describe the main equation, which follows from isoperimetry:
since at every step all neighbours of $K_n$ but the protected ones become burning, we have
\begin{equation}\label{eq:main} k_{n+1}\geq k_n + \Phi(k_n)-g(n) \end{equation}


\begin{proof}[\textbf{Proof of Theorem \ref{thm:main}}]
Any group $G$ satisfying  $v(n)\succeq n^d$ satisfies a $d$-dimensional isoperimetric inequality (see e.g. \cite{PSC}). That is $\Phi(k)\geq ck^{\frac{d-1}{d}}$ for some constant $c=c_G$. (In fact, if $v(n)\simeq n^d$ then $\Phi(k)\simeq  k^{\frac{d-1}{d}}$, but we will not need this).

The main equation \eqref{eq:main} now becomes
\begin{equation}
k_{n+1}\geq k_n + ck_n^{\frac{d-1}{d}} - g(n)
\end{equation}
with $c>0$ depending only on $\langle G,S\rangle$ and $g(n)=\sum_{k=1}^n f(k) = o(n^{d-1})$.

We now argue by induction that there exist $F,R$ depending only on $c,d$ and the function $g$ such that if $k_n>\frac{(n+F)^d}{R^d}$ then $k_{n+1}\geq \frac{(n+1+F)^d}{R^d}$.

To see this, first note that if $R>{4}{cd}$ then
\begin{multline*}
  \frac{(n+F+2)^d}{R^d} \leq \frac{(n+F)^d}{R^d}+\frac{2d (F+n)^{d-1}}{R^d}+O((F+n)^{d-2}) \\
  \leq \frac{(n+F)^d}{R^d}+\frac{c}{2} \frac{(F+n)^{d-1}}{R^{d-1}}+O((n+F)^{d-2}).
\end{multline*}
Therefore, if $F$ is large enough, $k_n + ck_n^{\frac{d-1}{d}} \geq \frac{(n+F+2)^d}{R^d}$.

Now by increasing $F$ further if needed we get  $$k_{n+1}\geq k_n + ck_n^{\frac{d-1}{d}} -g(n)\geq \frac{(n+F+2)^d}{R^d}-o(n^{d-1}) \geq \frac{(n+F+1)^d}{R^d}.$$
In particular, if the initial set of burning vertices $K_0$ satisfies $k_0>\frac{F^d}{R^d}$ then by induction $k_n\geq \frac{(n+F)^d}{R^d} \rightarrow \infty$.
Thus any big enough initial fire cannot be $\{f(n)\}$-contained.
\end{proof}


\begin{proof}[\textbf{Proof of Corollary \ref{c:inter}}]
By \cite{CSC} (See also \cite{PGG} Section 5), any group of intermediate growth satisfies a $d$-dimensional isoperimetric inequality $\Phi(k)\geq ck^{\frac{d-1}{d}}$ for all $d\in \NN$.
It follows form Theorem \ref{thm:main} that $G$ does not satisfy $\{f(n)\}$ containment for any $f=o(n^{d-2})$ for all $d$.
\end{proof}

\begin{proof} [\textbf{Proof of Theorem \ref{t:strectched}}]
By \cite{CSC,PGG}, $G$ satisfies the following isoperimetric inequality: $\Phi(k)\geq c\frac{k}{(\log k)^{1/\alpha}}$. So Equation \eqref{eq:main} becomes
$$k_{n+1}\geq k_n + c\frac{k_n}{(\log k_n)^{1/\alpha}} - g(n)$$
We now claim that if $\beta\geq \frac{\alpha}{\alpha+1}$, $f(n)=o(\frac{e^{n^\beta}}{n})$ and $k_n\geq Ce^{n^\beta}$ for some large enough $C$ (depending only on $c,f$) then $k_{n+1}\geq Ce^{(n+1)^\beta}$.
This would imply the theorem as it would mean that for any initial fire of size $k_0\geq C$,  $k_n\rightarrow \infty$.
To prove the inductive claim we may assume $n$ is large enough, use the fact that $e^{(n+2)^\beta}\leq e^{n^\beta}(1+\frac{3}{n^{1-\beta}})$ and note that the condition on $\beta$ implies $n^{1-\beta} > n^{\beta/\alpha}$. Therefore, for large enough $n$
$$k_n+c\frac{k_n}{(\log k_n)^{1/\alpha}}\geq Ce^{n^\beta}+c\frac{Ce^{n^\beta}}{2n^{\beta/\alpha}} \geq Ce^{(n+2)^\beta}.$$
Since $g(n)=\sum_{k=1}^n f(n)=o(e^{n^\beta})$ we conclude that $k_{n+1}\geq e^{(n+1)^\beta}$ as required.
The theorem follows as we can get rid of the $1/n$ factor in $f$ by slightly perturbing $\beta$.
\end{proof}

\begin{proof}[\textbf{Proof of Theorem \ref{t:strectched}}]
By \cite{EZ}, $G_\omega$ satisfies $\Phi(k)\geq \frac{ck}{\log k}$ for some $c>0$. Applying the proof of Theorem \ref{t:strectched} verbatim with $\alpha=1$ now gives the result.
\end{proof}

\subsection*{Acknowledgements}
While performing this research, G.A. and R.B. were supported by the Israel Science Foundation grant \#575/16 and by GIF grant \#I-1363-304.6/2016. G.K. was supported by the Israel Science Foundation grant \#1369/15 and by the Jesselson Foundation.

\end{document}